\documentclass[11pt,letter]{article}
\usepackage{amsfonts,amsmath,amssymb,graphicx,tabularx,booktabs}
\usepackage{mathtools}

\newtheorem{theorem}{\bf Theorem}[section]

\newtheorem{conjecture}[theorem]{\bf Conjecture}

\newcommand{\QED} {\hfill$\blacksquare$}

\title{\bf{On conjectures of network distance measures by
using graph spectra}}

\author{
Aleksandar Ili\'c \\
Facebook Inc, 1 Hacker Way, Menlo Park, 94025 California,USA\\
e-mail: \tt{aleksandari@gmail.com} 
\and
Matthias Dehmer \\
UMIT- The Health and Life Science University, \\
Department for Biomedical Computer Science and Mechatronics \\
6060 Hall in Tyrol, Austria \\
e-mail: \tt{matthias.dehmer@umit.at} 
}

\date{\today}

\begin{document}

\maketitle

\begin{abstract}
In this note we resolve three conjectures from [M. Dehmer, S. Pickl, Y. Shi, G. Yu,
\emph{New inequalities for network distance measures by using graph spectra}, Discrete Appl. Math. 252 (2019), 17--27] on the comparison of distance measures based on the graph
spectra, by constructing families of counterexamples and using computer search.
\end{abstract}

\vskip 0.1cm

 \vskip 0.1cm

\section{Introduction}
Eigenvalues of various graph-theoretical matrices often reflect structure properties of graphs
meaningfully \cite{cvetkovic1}. In fact, studying eigenvalues of graphs and, then, characterizing those graphs based on certain properties
of their eigenvalues has a long standing history, see \cite{cvetkovic1}. Also, eigenvalues have been used for characterzing graphs
quantitatively in terms of defining graph complexity as well similarity measures \cite{randic_2001_2,DeEm14}. An analysis revealed that
eigenvalue-based graphs measures tend to be quite unique, i.e., they are able to discriminate graphs uniquely 
\cite{dehmer_grabner_2012_2}. Some of the studied measures even outperformed measures from the family of the so-called Molecular ID Numbers, see
\cite{dehmer_grabner_2012_2}.

In this short paper, we further investigate an approach in \cite{DeEm14} where the authors explored inequalities for graph distance measures.
Those are
based on topological indices using eigenvalues of adjacency matrix, Laplacian matrix and signless
Laplacian matrix.
The graph distance measure is defined as
$$
d_I(G, H) = d(I(G), I(H) = 1 - e^{- \left( \frac{I(G) - I(H)}{\sigma} \right) ^2},
$$
where $G$ and $H$ are two graphs and $I(G)$ and $I(H)$ are topological indices applied to both $G$
and $H$.

In this short note, we are going to disprove three conjectures proposed in Dehmer et al. \cite{DePi19}, 
by constructing families of counterexamples and using computer search.

\section{Main result}

Let $G$ be a simple connected graph on $n$ vertices. Let $\lambda_1$ be the largest eigenvalue of the adjacency matrix of $G$, and $q_1$ be the largest eigenvalue of the Laplacian matrix of $G$.

The authors from \cite{DePi19} proposed the following conjectures and stated that it is likely that we need deeper results from matrix theory and from the theory of graph spectra to prove these.

\begin{conjecture}
Let $T$ and $T'$ be two trees on $n$ vertices. Then
$$
d_{q_1} (T, T')  \geq d_{\lambda_1} (T, T').
$$
\end{conjecture}

We are going to disprove the above conjecture by providing a family of counterexamples for
which it holds
$$
0 = d_{q_1} (T, T') < d_{\lambda_1} (T, T'),
$$
or in other words $q_1(T) = q_1(T')$ and $\lambda_1 (T) \neq \lambda_1 (T')$.

In \cite{Os13}, the author proved the following result: Almost all trees have a cospectral mate with
respect to the Laplacian matrix.

\begin{theorem}
Given fixed rooted graphs $(G, v)$ and $(H, v)$ and an arbitrary rooted graph $(K, w)$,
if $(G, u)$ and $(H, v)$ are Laplacian (signless Laplacian, normalized Laplacian, adjacency) cospectrally rooted then $G\cdot K$ and $H \cdot K$ are cospectral with respect to the Laplacian (signless Laplacian, normalized Laplacian, adjacency) matrix.
\end{theorem}

Starting from Laplacian cospectrally rooted trees shown in Figure 1 - one can construct
many graphs by choosing arbitrary trees $K$.

\begin{figure}[h]
\centering
\includegraphics[height=3cm]{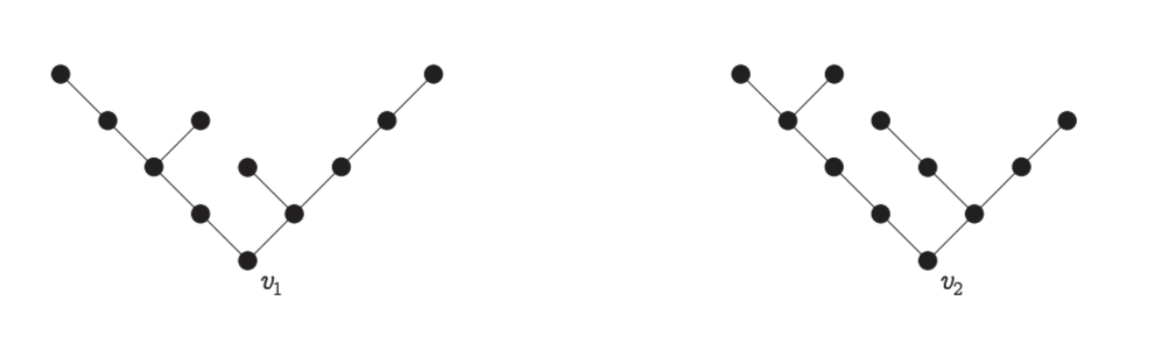}
\caption{Rooted Laplacian cospectral trees.}
\end{figure}

By direct calculation, we get that these trees are not adjacency cospectral and therefore the
adjacency spectral radiuses are different (2.0684 vs 2.0743). 
We rerun the same simulation for trees using \texttt{Nauty} \cite{Mc81} on $n = 10$ vertices 
as discussed in \cite{DePi19}. Based on the computed search - the smallest counterexample is on $n = 8$ vertices. 

\begin{center}
\begin{tabular}{ |r|r|r| } 
 \hline
 $n$ & tree pairs & Conjecture 5.1 counterexamples \\ 
\hline 
4 & 3 & 0 \\
5 & 6 & 0 \\
6 & 21 & 0 \\
7 & 66 & 0 \\
8 & 276 & 2 \\
9 & 1128 & 11 \\
10 & 5671 & 89 \\
11 & 27730 & 568 \\
12 & 152076 & 3532 \\
13 & 846951 & 21726 \\
14 & 4991220 & 138080 \\
15 & 29965411 & 877546 \\
16 & 186640860 & 5725833 \\
 \hline
\end{tabular}
\end{center}

Degree powers, or the zeroth Randi\'{c} index are defined as
$$
F_k = \sum_{v \in V} deg^k (v).
$$

\begin{conjecture}
Let $T$ and $T'$ be two trees on $n$ vertices. Then
$$
d_{F_2} (T, T')  \geq d_{q_1} (T, T').
$$
\end{conjecture}

We rerun the same computer simulation and found many examples of pairs for which holds 
$$
|F_2 (T) - F_2(T')|  < |q_1(T) - q_1(T')|,
$$
and consequently $d_{F_2} (T, T')  < d_{q_1} (T, T')$. In particular the smallest counterexample is on $n=6$ vertices and shown in Figure 2: clearly $F_2(T) = F_2(T') = 20$ and 
$$4.214320 = q_1(T) < q_1(T') = 4.302776.$$ 

This disproves the above conjecture and corrects the results from \cite{DePi19}.

\begin{center}
\begin{tabular}{ |r|r|r| } 
 \hline
 $n$ & tree pairs & Conjecture 5.2 counterexamples \\ 
\hline 
4 & 3 & 0 \\
5 & 6 & 0 \\
6 & 21 & 1 \\
7 & 66 & 5 \\
8 & 276 & 28 \\
9 & 1128 & 117 \\
10 & 5671 & 577 \\
11 & 27730 & 2672 \\
12 & 152076 & 13805 \\
13 & 846951 & 72801 \\
14 & 4991220 & 405454 \\
15 & 29965411 & 2312368 \\
16 & 186640860 & 13713949 \\
 \hline
\end{tabular}
\end{center}

\begin{figure}[h]
\centering
\includegraphics[height=3cm]{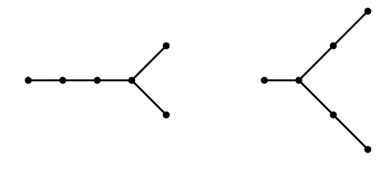}
\caption{Two trees with $F_2(T) = F_2(T')$ and $q_1(T) \neq q_1 (T')$.}
\end{figure}

To conclude, these results disprove Conjectures 5.1 and Conjecture 5.2 from \cite{DePi19}, 
while Conjecture 5.3 on relationship between $\lambda_1$ and $F_2$ \cite{DePi19} directly follows from these.

\section{Acknowledgement}
Matthias Dehmer thanks the Austria Science Funds (FWF) for financial support (P 30031).


\begin{thebibliography}{99}

\bibitem{cvetkovic1} D.~M. Cvetkovi\'c, M. Doob, H. Sachs, Spectra of Graphs. Theory and Application,
Deutscher Verlag der Wissenschaften, 1980, Berlin, Germany

\bibitem{DePi19}
M. Dehmer, S. Pickl, Y. Shi, G. Yu, \textit{New inequalities for network distance measures by
using graph spectra}, Discrete Appl. Math. 252 (2019), 17--27.

\bibitem{DeEm14}
M. Dehmer, F. Emmert-Streib, Y. Shi, \textit{Interrelations of graph distance measures based on
topological indices}, PLoS ONE, 9 (2014), e94985.

 
 \bibitem{dehmer_grabner_2012_2} M. Dehmer, M. Grabner, The Discrimination Power of Molecular Identification Numbers Revisited,
  MATCH Commun. Math. Comput. Chem., Vol. 69 (3), 2013, 785-794.
  
  \bibitem{Mc81}
B. D. McKay, \textit{Practical graph isomorphism}, Eur. J. Oper. Res. 30 (1981), 45--87.


\bibitem{randic_2001_2} M. Randi\'{c}, M. Vracko, M. Novi\' c, 
 Eigenvalues as Molecular Descriptors, In: QSPR/QSAR Studies by Molecular Descriptors, Editor: M. V. Diudea, 
 Nova Publishing, 2001, Huntington, NY, USA, 93--120.


\bibitem{Os13}
S. Osborne, \textit{Cospectral bipartite graphs for the normalized Laplacian}, Ph.D. Dissertation,
Iowa State University, Ames, 2013.


\end{thebibliography}
\end{document}